\documentclass[11pt]{amsart}
\usepackage{amsmath}
\usepackage{amssymb}
\usepackage{amscd}
\usepackage[dvips]{epsfig}
\usepackage{graphics}

\theoremstyle{plain}
\newtheorem{theorem}{Theorem}[section]
\newtheorem{lemma}[theorem]{Lemma}

\newtheorem{question}[theorem]{Question}

\theoremstyle{definition}
\newtheorem{definition}[theorem]{Definition}
\newtheorem{example}[theorem]{Example}

\newtheorem{remark}[theorem]{Remark}
\theoremstyle{remark}

\begin{document}
\def \a {\alpha}
\def \s {\sigma}

\title{Invariants of Legendrian Knots from Open Book Decompositions} 
\author{S\.{i}nem \c{C}el\.{i}k Onaran}

 \address{Department of Mathematics, Middle East Technical University, Ankara, Turkey and School of Mathematics, Georgia Institute of Technology, Atlanta, Georgia. } 
 \email{sinemo@math.gatech.edu, e114485@metu.edu.tr}
  \subjclass{57R17}
\keywords{contact structures, Legendrian knots, open book decompositions}
 \begin{abstract}In this note, we define a new invariant of a Legendrian knot in a contact $3$-manifold using an open book decomposition supporting the contact structure. We define the support genus $sg(L)$ of a Legendrian knot $L$ in a contact $3$-manifold $(M, \xi)$ as the minimal genus of a page of an open book of $M$ supporting the contact structure $\xi$ such that $L$ sits on a page and the framings given by the contact structure and the page agree. We show any null-homologous loose knot in an overtwisted contact structure has support genus zero. To prove this, we show that any topological knot or link in any $3$-manifold $M$ sits on a page of a planar open book decomposition of $M$.
 \end{abstract}
 \maketitle
 
 \setcounter{section}{0}
  \section{Introduction} 
\par Recently, contact geometry has been a major development in low dimensional topology due to work of Eliashberg, Giroux, Etnyre, Honda and many other mathematicians. The study of Legendrian knots is important in the theory since Legendrian knots reveal the geometry and topology of the underlying contact $3$-manifold. For example, Legendrian knots are used to distinguish contact structures \cite{K}, to detect topological properties of knots \cite{R} and to detect overtwistedness of contact structures \cite{EH2}.
\par In \cite{JB}, given any $3$-manifold $M$, Etnyre and Ozbagci defined new invariants of contact structures on $M$ in terms of open book decompositions supporting the contact structure. One of the invariants is the support genus of the contact structure which is defined as the minimal genus of a page of an open book of $M$ that supports the contact structure. In a similar fashion, we define the support genus $sg(L)$ of a Legendrian knot $L$ in a contact $3$-manifold $(M, \xi)$ as the minimal genus of a page of an open book of $M$ supporting the contact structure $\xi$ such that $L$ sits on a page and the framings given by the contact structure and the page agree. This definition is originally due to Etnyre.
\par First we study the topological properties of links sitting on pages of open book decompositions. For links in $S^3$, we explicitly construct a planar open book decomposition of $S^3$ which contains the link on its page.
\begin{theorem} Any topological link in $S^3$ sits on a planar page of an open book decomposition of $S^3$ whose monodromy is a product of positive Dehn twists. 
\end{theorem}
As a consequence of this theorem, we have a general property for topological links, in particular for knots. 
\begin{theorem}Any topological link in a closed, orientable $3$-manifold $M$ sits on a planar page of an open book decomposition of $M$.\label{manifold}
\end{theorem}
\par Next we study the following question: What contact geometric properties of knots are reflected by topological properties of knots sitting on pages of open books? Using the above theorems, we prove:
\begin{theorem} Any null-homologous loose Legendrian knot in an overtwisted contact $3$-manifold has support genus $sg(L) = 0$. \label{loose}
\end{theorem}
\par In the following section, we give a review of background information on contact structures, Legendrian knots in contact manifolds and open book decompositions. In Section 3, for a given topological link in $S^3$ we present an explicit algorithm to construct a planar open book decomposition which contains the given link on its page. We prove Theorem ~ \ref{manifold}. Using this, in Section 4, we conclude the proof of Theorem ~ \ref{loose}. Finally, in Section 5, we list several observations and open problems related to the support genus of Legendrian knots in contact $3$-manifolds. We show that we may arrange the monodromy of a planar open book decomposition for links in $S^3$ to be a product of positive Dehn twists. We conclude that for any given knot type $K$ in $(S^3, \xi_{std})$, there is a Legendrian representative $L$ of $K$ such that $sg(L) = 0$. We construct examples of non-loose Legendrian knots having support genus zero or non-zero. We show the existence of Legendrian knots with non-zero support genus in weakly fillable contact structures. Moreover, we observe that for a non-zero rational number $r \in \mathbb{Q}$, any contact $3$-manifold which is obtained by a contact $r$-surgery on a support genus zero Legendrian knot has support genus zero. 

\section{Background Information} 
Let us now review the basics of contact geometry and briefly mention the facts that we will use throughout the paper. 
\subsection{Contact Structures and Knots in Contact Manifolds} A \emph{contact structure} $\xi$ on an oriented $3$-manifold $M$ is a maximally non-integrable $2$-plane field. Locally $\xi$ can be given as a kernel of a $1$-form $\a$ and from the non-integrability condition we have $\a\wedge d\a \neq 0$. If $\xi$ is orientable, in this case $1$-form $\alpha$ exists globally and the $1$-form $\a$ is called a \emph{contact form}. We denote a \emph{contact $3$-manifold} as ($M,\xi$). 
\par  A knot $L$ in a contact $3$-manifold ($M,\xi$) is called \emph{Legendrian} if it is everywhere tangent to $\xi$. The classical invariants of Legendrian knots are the topological knot type, the Thurston-Bennequin invariant $tb(L)$ and the rotation number$rot(L)$. \emph{The Thurston-Bennequin invariant} $tb(L)$ measures the framing of $L$ given by the contact planes with respect to the framing given by a Seifert surface of $L$. \emph{The rotation number} $rot(L)$ of an oriented null-homologous Legendrian knot $L$ can be computed as the winding number of $TL$ after trivializing $\xi$ along a Seifert surface for $L$. 
\par
\begin{figure}[h!]
\begin{center}
  \includegraphics[width=8cm]{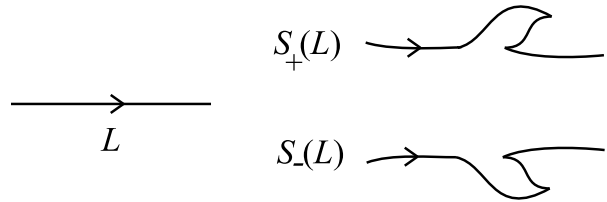}
 \caption{The positive stabilization $S_{+}(L)$ and the negative stabilization $S_{-}(L)$ of $L$.}
  \label{stab}
\end{center}
\end{figure}

 \emph{Positive (negative) stabilization} $S_{+}(L)$ ($S_{-}(L)$) of a Legendrian knot $L$ in the standard contact structure $\xi_{std}$ on $\mathbb{R}^3$ is obtained by modifying the front projection of $L$ by adding a down cusp (an up cusp) as in Figure ~ \ref{stab}. Since stabilizations are done locally, by Darboux's theorem this defines stabilizations of Legendrian knots in any contact $3$-manifold ($M,\xi$). After stabilizing a Legendrian knot the classical invariants change as $tb(S_{\pm}(L)) = tb(L)- 1$ and $rot(S_{\pm}(L)) = rot(L) \pm 1$.

\par A contact structure $\xi$ on $M$ is \emph{overtwisted} if there is an embedded disk with a Legendrian boundary having Thurston-Bennequin invariant zero, otherwise $\xi$ is called \emph{tight}. A Legendrian knot in an overtwisted contact $3$-manifold $M$ is \emph{loose} if its complement is also overtwisted. We call a Legendrian knot \emph{non-loose} if its complement is tight.
\par For details on contact structures and knots in contact structures, see \cite{J1}, \cite{J2}. 

\subsection{Open Book Decompositions}
 An \emph{open book decomposition} of a closed, oriented $3$-manifold $M$ is a pair ($S, \varphi$) where $S$ is an oriented compact surface with boundary link $B$ and \emph{the monodromy} map $\varphi$ is a diffeomorphism of $S$ such that $\varphi$ is identity on a neighborhood of the boundary $\partial S$ and $M - B = S \times [0, 1] / (1, x) \sim (0, \varphi(x))$. The fibers are the interior of the Seifert surface $S$ of $B$. The link $B$ is called the \emph{binding} and the fiber surface $S$ is called the \emph{page} of the open book decomposition. The \emph{genus} of an open book decomposition is defined as the genus of the page. In particular, planar open book decompositions are genus zero open book decompositions.

\par \emph{Positive (negative) stabilization} of an open book decomposition ($S, \varphi$) is the open book decomposition ($S', \varphi\circ{t_a}^{\pm 1}$) where $S' = S \cup (1$-handle) and $t_a$ ($t_a^{-1}$) is a right (left) handed Dehn twist along the closed curve $a$ in $S'$ running over the $1$-handle and intersecting the co-core of the $1$-handle once.

\par An open book decomposition of $M$ and a contact structure $\xi$ on $M$ are \emph{compatible} if after an isotopy of the contact structure, there is a contact form $\a$ for $\xi$ such that $\a > 0$ on the binding $B$, in other words the binding $B$ is a positive transverse link, and $d\a > 0$ on every page of the open book decomposition. 
\par In \cite{TW}, Thurston and Winkelnkemper show that every open book decomposition of a $3$-manifold admits a compatible contact structure. In \cite{G}, Giroux proves that every contact structure is compatible with some open book decomposition and he also proves that there is a one to one correspondence between oriented contact structures up to isotopy and open book decompositions up to positive stabilization.
\par
\begin{figure}[h!]
\begin{center}
  \includegraphics[width=9cm]{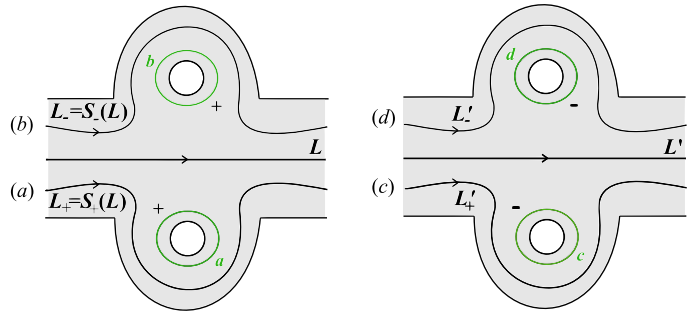}
 \caption{$(a)$ Positive stabilization $S_{+}(L)$ of $L$, $(b)$ Negative stabilization $S_{-}(L)$ of $L$, $(c)$ Positive destabilization of $L'$, $S_{+}(L'_{+}) = L'$, $(d)$ Negative destabilization of $L'$, $S_{-}(L'_{-})= L'$.}
  \label{obstab}
\end{center}
\end{figure}

 The following lemma is useful and gives the relation between the stabilizations of open book decompositions and the stabilizations of Legendrian knots sitting on a page of an open book decomposition. For the proof see Appendix.
\begin{lemma}Let ($S, \varphi$) be an open book decomposition for a closed oriented $3$-manifold $M$ compatible with a contact structure $\xi$ on $M$. Let $L$ be a Legendrian knot sitting on a page of the open book.
\begin{enumerate}
\item Positive (resp. negative) stabilization $S_{+}(L)$ (resp. $S_{-}(L)$) of the Legendrian knot $L$ can be realized on the page of the open book by first stabilizing the open book positively and then pushing the knot $L$ over the $1$-handle that we use to stabilize the open book. See Figure ~ \ref{obstab}$(a)$ and $(b)$.
\item If we first negatively stabilize the open book and then push the knot $L$ over the $1$-handle that we use to stabilize the open book, then the negatively stabilized open book is no longer compatible with the contact structure $\xi$, but the curve $L$ on the page gives a Legendrian knot $L'$ in the new contact structure and Legendrian knots $L'_{+}$ and $L'_{-}$ in Figure~\ref{obstab}$(c)$ and $(d)$ are positive and negative destabilizations of $L'$, respectively.
\end{enumerate}
\label{rkstab}
\end{lemma}
\par To prove the main theorems we also need the following lemma.
\begin{lemma} Let $M$ be a closed oriented $3$-manifold and let $(S, \varphi)$ be an open book decomposition for $M$.
\begin{enumerate} 
\item If $K$ is a knot in $M$ intersecting each page $S$ transversely once, then the result of a $0$-surgery along $K$ gives a new manifold with an open book decomposition having a page $S' = S -  $\{open disk\} and having the knot $K$ as one of the binding components. In particular, if the knot $K = \{x\} \times [0, 1] /  \sim $ in the mapping torus $M_{\varphi}$ in $M$ for a fixed point $x \in S$ of $\varphi$ and if $\varphi \mid_{\{open\,disk\}} = id$ then the new monodromy $\varphi'$ after a $0$-surgery along K is $\varphi' = \varphi \mid_{{\tiny S'}}$. 
\item If $K$ is a knot in $M$ sitting on a page $S$ of the open book decomposition, then $\pm1$-surgery along $K$ with respect to the page framing gives a new manifold with an open book decomposition $(S, \varphi \circ t^{{\tiny \mp1}}_{K})$ where $t^{{\tiny +1}}_{K}$/ $t^{{\tiny -1}}_{K}$ denotes right/ left handed Dehn twists along the knot $K$.
\end{enumerate}
\label{lemma2}
\end{lemma}
\par For the proof of above Lemma ~ \ref{lemma2} and for more information on open book decompositions, see \cite{J3}.
\section{Invariants of Knots from Open Book Decompositions}
\begin{theorem} Any knot $K$ in $S^3$ is planar, that is $K$ sits on a page of a planar open book decomposition for $S^3$.
   \label{top}
 \end{theorem}
 Before the proof of Theorem ~ \ref{top}, let us define some terminology, give an illustrative example and state a fundamental lemma.
\par
\begin{figure}[h!]
\begin{center}
  \includegraphics[width=10cm]{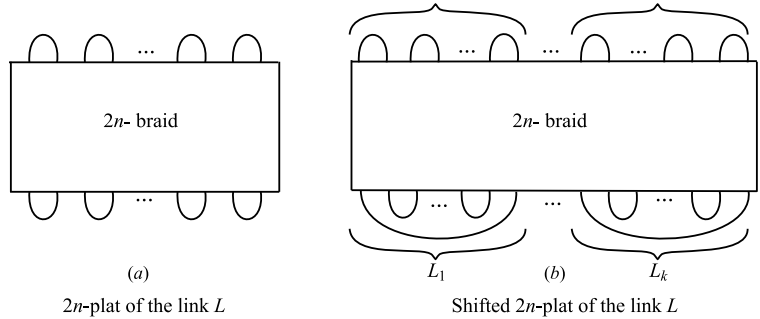}
 \caption{Shifted $2n$-plat of $L$ is called \emph{pure braided $2n$-plat} if its associated $2n$-braid is a pure braid.}
  \label{plat}
\end{center}
\end{figure}

\begin{figure}[h!]
\begin{center}
  \includegraphics[width=13cm]{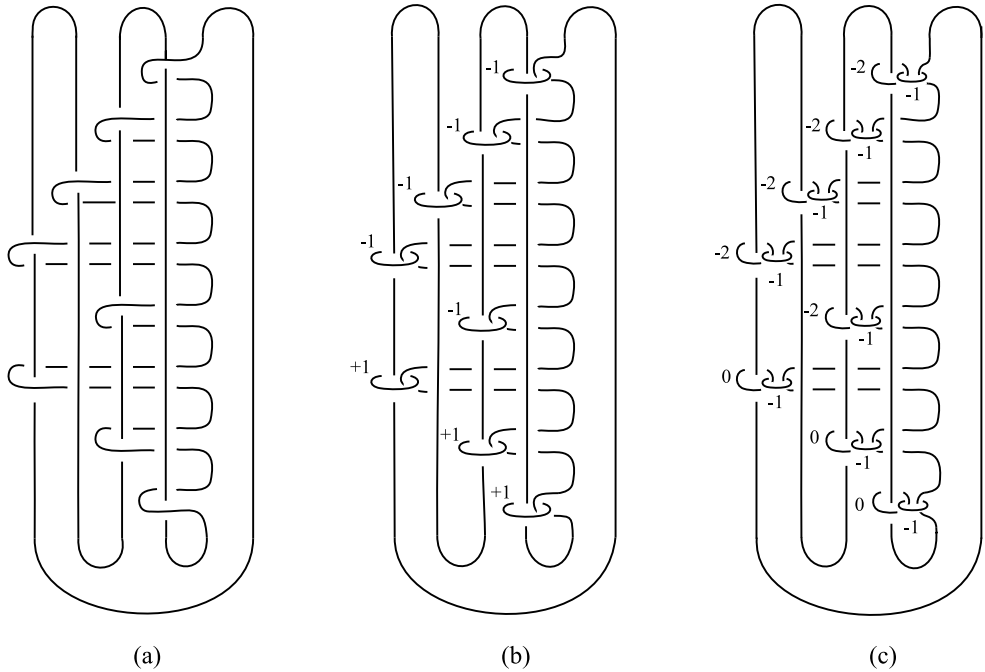}
 \caption{Unknotting the figure eight knot.}
  \label{puregen}
\end{center}
\end{figure}

 It is well known that any link $L$ of $k$ components $L_1, \ldots, L_k$, in particular any knot $K$, can be represented as a $2n$-plat, see Figure ~ \ref{plat}$(a)$. We define the \emph{shifted $2n$-plat} of the link $L$ as the closure of a $2n$-braid as shown in Figure ~ \ref{plat}$(b)$. We say a shifted $2n$-plat of the link $L$ is \emph{pure braided $2n$-plat} if its associated $2n$-braid is a pure braid. For details on braid group, pure braid group and plat presentation of knots and links, see \cite{B}.
\begin{example} The figure eight knot $K$ is planar. The aim here is to present the figure eight knot $K$ as a pure braided plat as in Figure ~ \ref{puregen}$(a)$ and using this pure braided plat and the ideas in Lemma ~ \ref{lemma2} to construct a planar open book which contains the figure eight knot on its page. 
  \begin{figure}[h!]
\begin{center}
  \includegraphics[width=2.4cm]{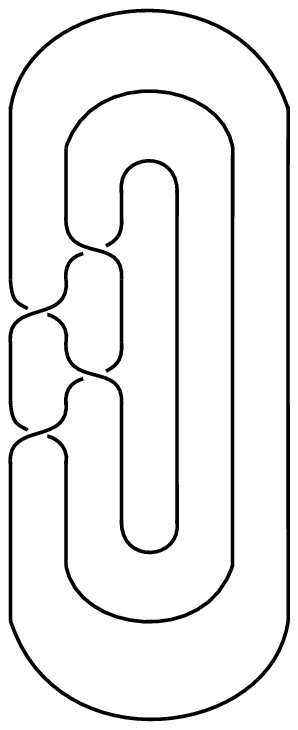}
 \caption{Braid representative of the figure eight knot.}
  \label{f8}
\end{center}
\end{figure}
\end{example} 
\par
\begin{figure}[h!]
\begin{center}
  \includegraphics[width=13cm]{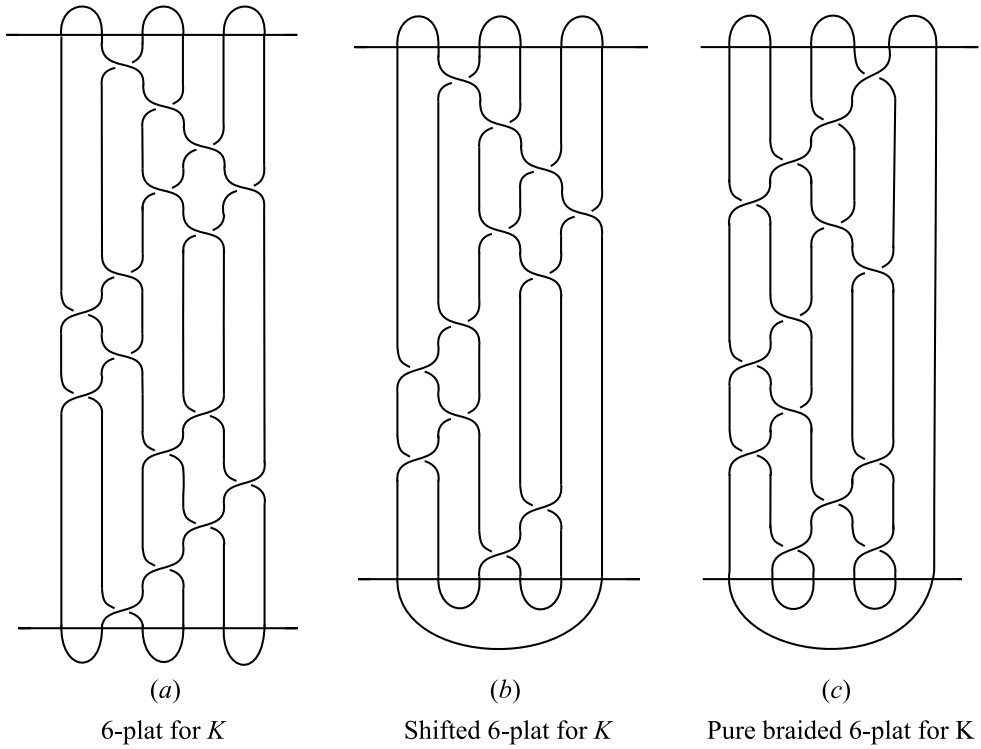}
 \caption{Pure braided plat presentation of the figure eight knot.}
  \label{pure}
\end{center}
\end{figure}

 We start with a minimum braid representation of the figure eight knot $K$ as in Figure ~ \ref{f8}. Throughout $\s_i$, $i = 1, \ldots, n-1$, stand for the standard generators of the braid group $B_n$ on $n$-strands. Note that $K$ has braid index $3$ and its associated braid word is $b = \s_2^{-1}\s_1\s_2^{-1}\s_1$. As shown in Figure ~ \ref{pure}$(a)$, we can represent $K$ by a $6$-plat which is associated to a $6$-braid $b_0\tilde{b}{b_0}^{-1}$ where $b_0 = (\s_2\s_3\s_4\s_5)(\s_3\s_4)$ and $\tilde{b}$ is the $6$-braid obtained from $b$ by adding $3$ trivially braided strands. Now isotope the diagram in Figure ~ \ref{pure}$(a)$ to obtain a shifted $6$-plat as in Figure ~ \ref{pure}$(b)$ and continue isotoping to obtain a pure braided $6$-plat for the figure eight knot as in Figure ~ \ref{pure}$(c)$. In Lemma ~ \ref{lemma1} below, we present an algorithm to obtain a pure braided plat for any given knot or link in $S^3$.
\par
\begin{figure}[h!]
\begin{center}
  \includegraphics[width=11cm]{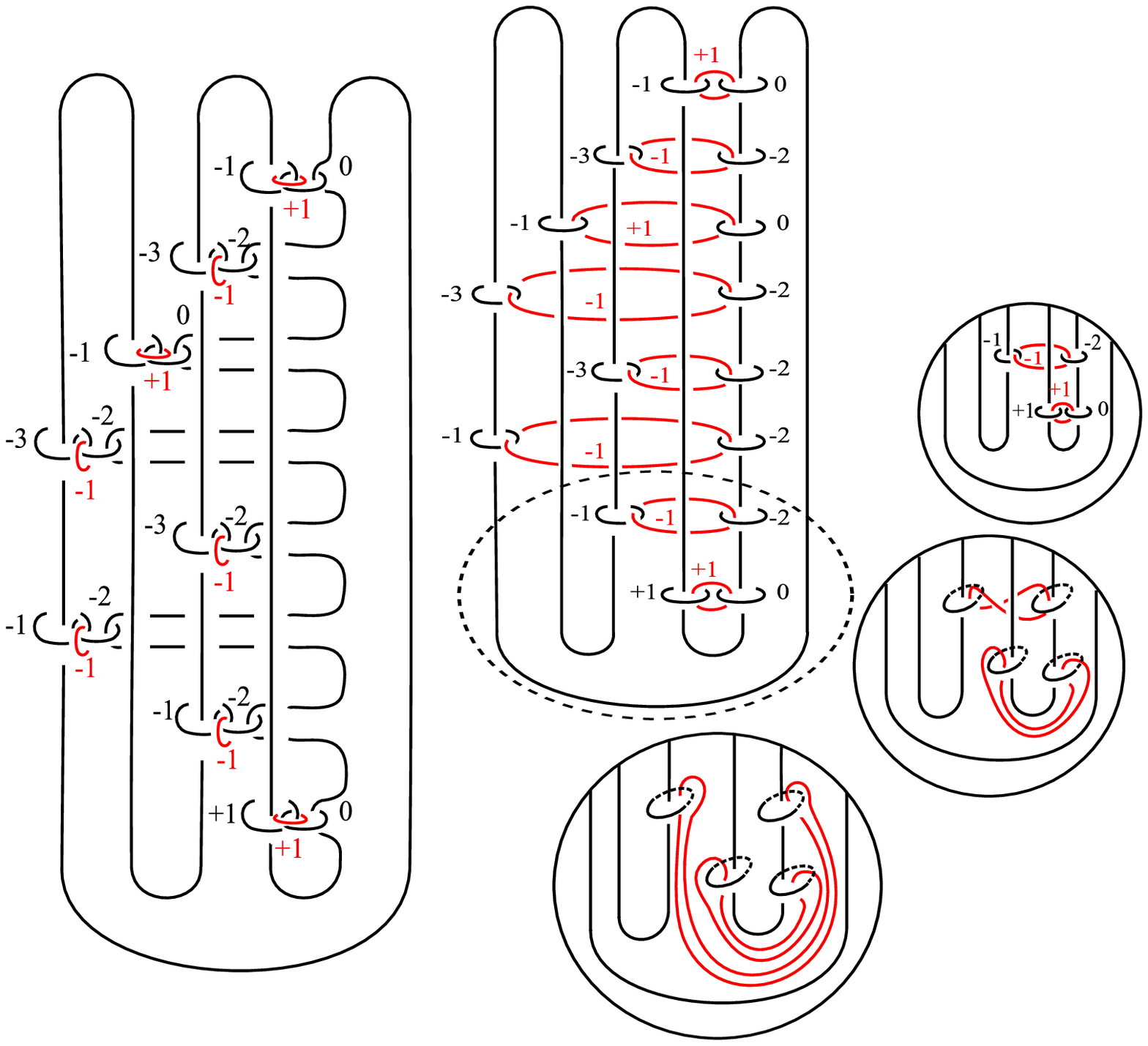}
 \caption{The unknotted knot $K$ bounds a disk and we isotope the middle $\pm1$-framed unknots onto the disk.}
  \label{3}
\end{center}
\end{figure}

 Next, we decompose the pure braided $6$-plat of the figure eight knot in standard generators of the pure braid group on $6$-strands as in Figure ~ \ref{puregen}$(a)$. Now to obtain the open book decomposition which contains the figure eight knot $K$, we unknot $K$ using the diagram in Figure ~ \ref{puregen}$(a)$. We unknot $K$ by blowing up twists. See Figure ~ \ref{puregen}$(b)$. We get a link $L_K$ of unknots linking $K$ whose components have framing $\pm1$. We continue blowing up to ensure that each component of $L_K$ links $K$ exactly once. See Figure ~ \ref{puregen}$(c)$. Notice that we add new $\pm 1$-framed components to the link $L_K$ and the components of $L_K$ link each other as the Hopf link and link the knot $K$ only once. We continue blowing up as in Figure ~ \ref{3} to remove each linking between the components. We need to be careful with the resulting $\pm 1$-framed unknots linking the components of $L_K$. To be more precise, at each linking crossing between the components of $L_K$ we have different choices where to blow up as explained in the proof of Theorem ~ \ref{top} below. We always choose the one that guarantees that after blowing up, the resulting $\pm 1$-framed unknots linking the components of $L_K$ can be isotoped to sit on the page of the open book decomposition at the end. See Figure ~ \ref{3} again. 
\par
\begin{figure}[h!]
\begin{center}
  \includegraphics[width=11.5cm]{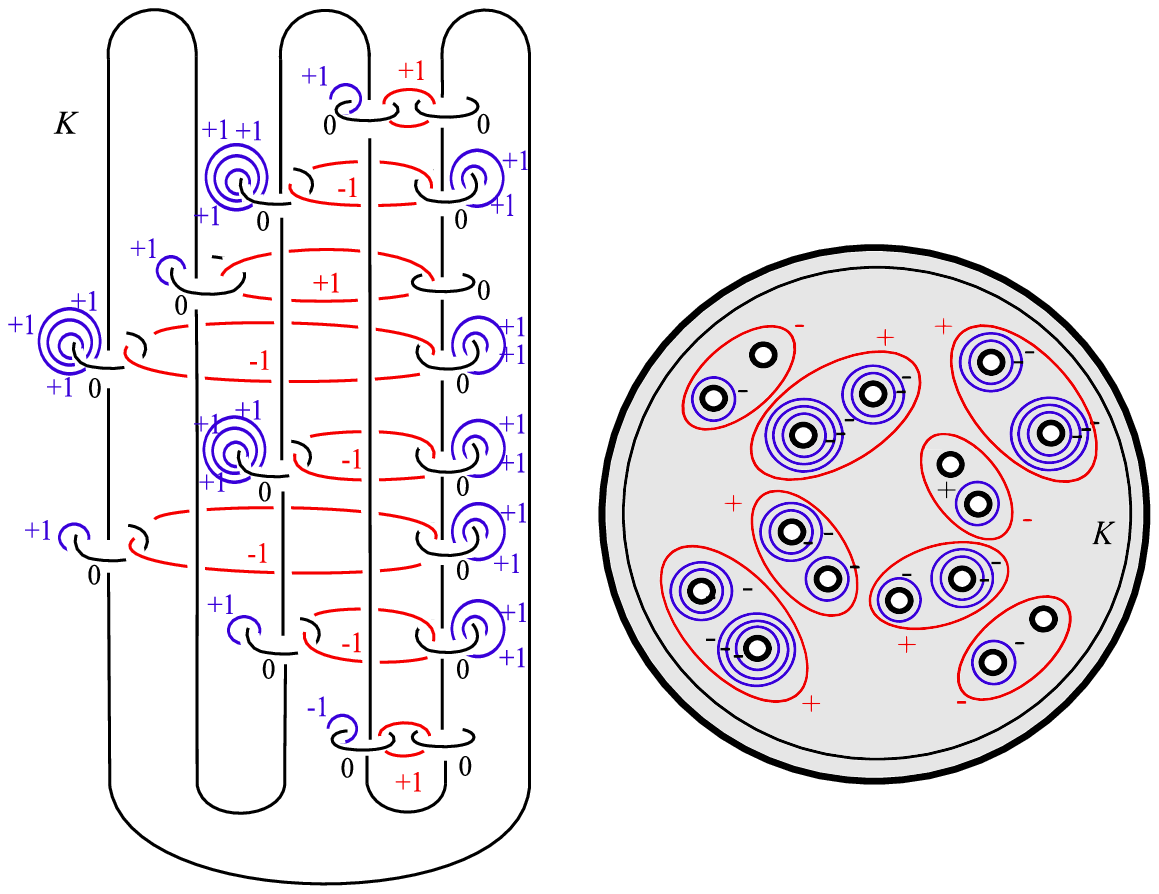}
 \caption{Page of a planar open book decomposition containing the figure eight knot, pages are disk with $16$ punctures.}
  \label{final}
\end{center}
\end{figure}

 Finally, we blow up again as in Figure ~ \ref{final} so that each component of the link $L_K$ has framing coefficient $0$. Now, using Lemma ~ \ref{lemma2} we are in a position to see the open book decomposition explicitly. Note that after performing surgeries, we obtain a planar open book decomposition for $S^3$ where the figure eight knot $K$ and each $0$-framed components of $L_K$ are the binding components of the open book decomposition. Each $\pm 1$-framed unknots linking the components of $L_K$ sits on the page and contributes negative/ positive Dehn twists to the monodromy of the open book decomposition respectively.
 \begin{lemma} $(1)$ Every knot can be represented as a pure braided plat.
 \smallskip
\par $(2)$ Every link of $k$ components $L_1, \ldots, L_k$ can be represented as a pure braided plat.
  \label{lemma1}
 \end{lemma}
   \begin{proof} $(1)$ We may isotope a shifted $2n$-plat of the knot $K$ to get a pure braided $2n$-plat for $K$ as follows: First orient the knot $K$ and label the lower and the upper end points of the strands of associated $2n$-braid $b$ and pair them as in Figure ~ \ref{ispat}. We have the following list of pairs: for the lower end points ($2n, 1$), ($2, 3$), \ldots, ($2n-2, 2n-1$) and for the upper end points ($1', 2'$), \ldots, ($(2n-1)', (2n)'$). Also, denote the permutation in the permutation group $S_{2n}$ on the set $\{1, \ldots, 2n\}$ associated to $2n$-braid $b$ of the shifted $2n$-plat by $\s$. 
\begin{figure}[hbt]
\begin{center}
  \includegraphics[width=9cm]{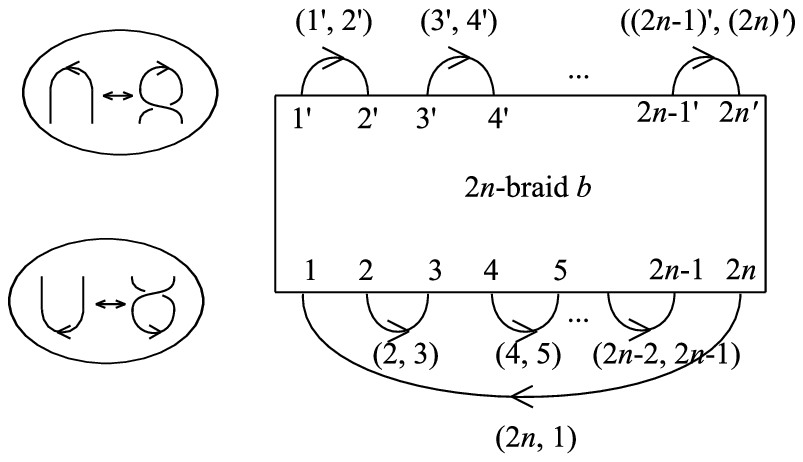}
 \caption{From a shifted $2n$-plat to a pure braided plat.}
  \label{ispat}
\end{center}
\end{figure} 
  \begin{figure}[hbt]
\begin{center}
  \includegraphics[width=14cm]{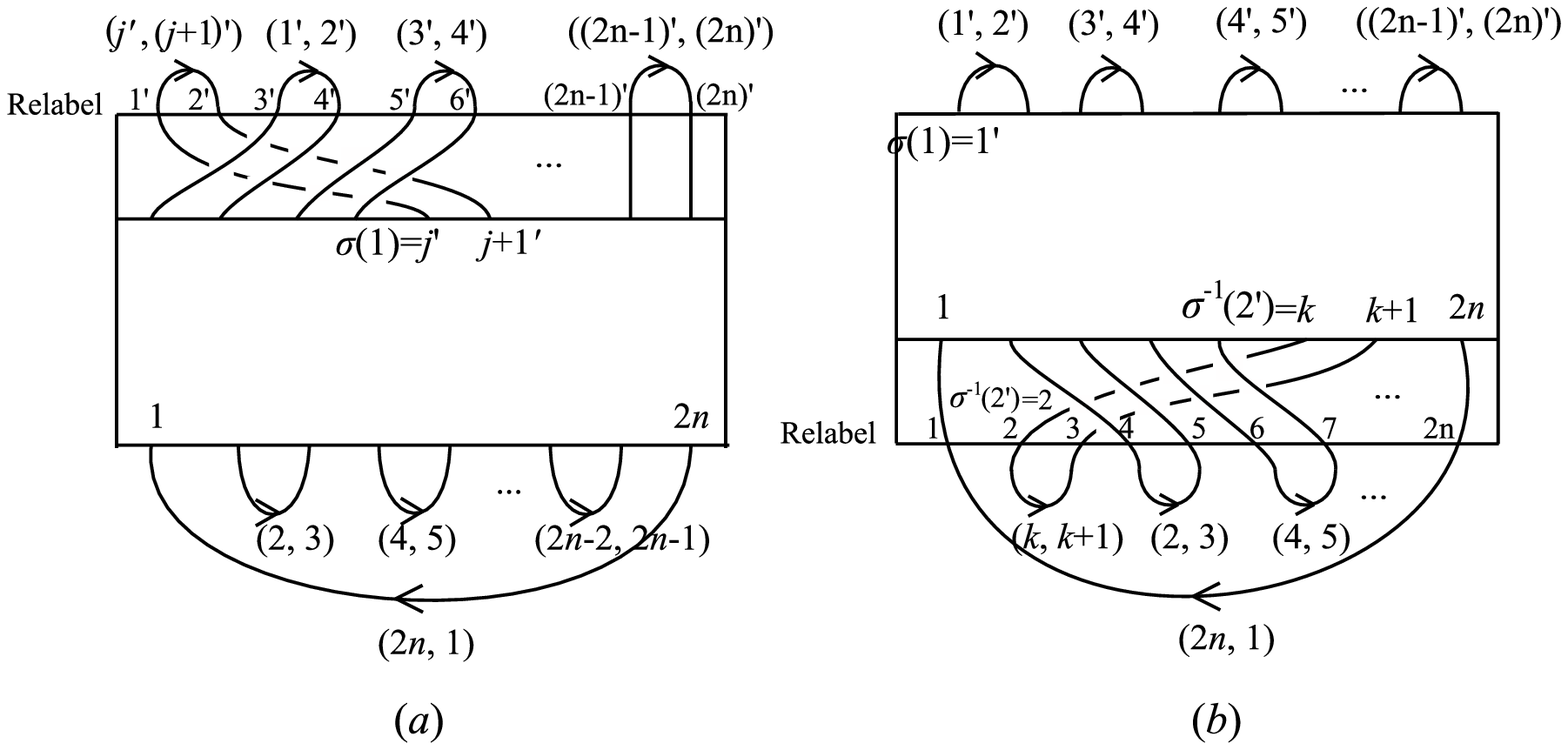}
 \caption{From a shifted $2n$-plat to a pure braided plat.}
  \label{ispat2}
\end{center}
\end{figure} 
     \par Now, start in the lower left strand with a labeled $1$ lower end point. This strand connects to its upper point $j' = \s(1)$. Isotope $(j', (j+1)')$ to the left as in Figure ~ \ref{ispat2}$(a)$ so that the first labeled upper point at the top is $j' = \s(1)$. Now relabel upper end points as $1', 2', \ldots, (2n)'$ and without loss of generality denote the permutation associated to new $2n$-braid as $\s$ again. Next, find where the strand whose upper end point is $2'$ connects at the bottom, its lower end point will be $\s^{-1}(2') = k$ where $k$ is an element from the set $\{2, \ldots, 2n-1\}$. Note that $k\neq2n$, otherwise the knot $K$ would be a link. Isotope $(k, k+1)$  to the left as in Figure ~ \ref{ispat2}$(b)$ to be the second labeled strand at the bottom. Relabel the lower end points as $1, 2, \ldots, 2n$ and without loss of generality again denote the permutation associated to new $2n$-braid as $\s$. Note that we have $\s(1) = 1'$, $\s(2') = 2$. Find $\s(3)$ and isotope similarly to be the third labeled strand at the top. Continuing in this manner, we will obtain a pure braid giving a pure braided $2n$-plat of the knot $K$.
\par $(2)$ First of all, given a link $L$ of $k$ components $L_1, \ldots, L_k$ we can present the link $L$ as a plat. From this plat we can obtain a shifted $2n$-plat of $L$ such that it has the same form as in Figure ~ \ref{plat}$(b)$ with an associated braid $b$ which is not necessarily a pure braid. However, the algorithm described in proof of $(1)$ extends to convert a shifted $2n$-plat of $L$ into a pure braided $2n$-plat.
 \end{proof}
{\em \textbf{Proof of Theorem ~ \ref{top}.}} Given a knot $K$ in $S^3$, we construct a planar open book of $S^3$ such that $K$ is one of the binding components. We then push the knot $K$ onto one of the pages.
\par
\begin{figure}[hbt]
\begin{center}
  \includegraphics[width=12cm]{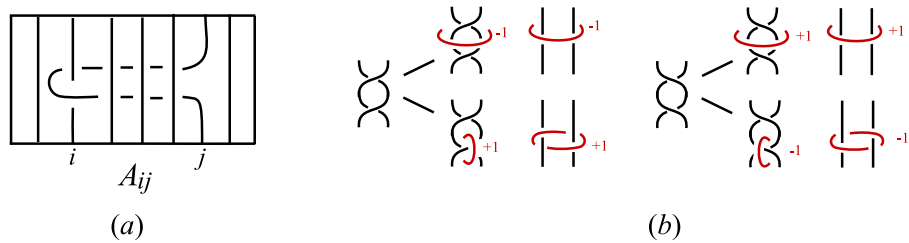}
 \caption{$(a)$ Generator $A_{ij}$ for the pure braid group, $(b)$ Different ways of blowing up to remove twists.}
  \label{unknotting}
\end{center}
\end{figure}

 First, present the knot $K$ as a pure braided plat using the algorithm given in Lemma ~ \ref{lemma1}. Next, decompose the pure braided plat of $K$ in terms of standard generators of the pure braid group. A generating set of braids $A_{ij}$, $ 1 \leq i < j \leq 2n$, for the pure braid group on $2n$-strands is shown in Figure \ref{unknotting}$(a)$. Note that to unknot the knot $K$ using a decomposed pure braided plat presentation of $K$, we only need to remove full twists. We remove twists and unknot $K$ by blowing up. Note also that there is not a unique way to do so. The different ways of blowing up are shown in Figure ~ \ref{unknotting}$(b)$. 
\par Idea of the proof is that using this presentation of the knot $K$, unknot $K$ by blowing up several times in such a way that $K$ is the unknot which we denote by $U_K$ and the resulting link of unknots $L_K$ coming from the blow ups step by step satisfy:
\begin{enumerate}
\item Each component of $L_{K}$ links $U_K$ only once,
\item The components of $L_{K}$ are pairwise unlinked or linked as the Hopf link, 
\item If the components of $L_K$ linked as the Hopf link, then continue blowing up to remove the linking and get $\pm1$-framed unknots $L_{\pm}$ linking the components of $L_K$,
\item $L_{\pm}$ does not link $U_K$ and each can be isotoped to sit on a disk that $U_K$ bounds,
\item The component of $L_{K}$ linking $U_K$ only once has $0$-framing, we denote such components by $L_0$.
\end{enumerate}
\par Thus, at the end we have a union of links of unknots $L_K = L_{0} \cup L_{\pm}$. Note that the unknot $U_K$ has a natural open book in $S^3$ coming from the disk it bounds. The $0$-framed link $L_{0}$ of unknots puncture each disk page transversely once and we can isotope $\pm 1$-framed link $L_{\pm}$ of unknots linking $L_0$ components of $L_K$ onto one of the punctured disk pages. Thus, after performing surgeries $U_K$ will be isotopic to the knot $K$ and by Lemma ~ \ref{lemma2} we will get a planar open book of $S^3$ where the knot $U_K$ and the $0$-framed link $L_{0}$ of unknots form the binding components and each $\pm 1$-framed link $L_{\pm}$ of unknots sitting on the punctured disk page contributes to negative/ positive Dehn twist to the monodromy of the new open book.
\par
\begin{figure}[h!]
\begin{center}
  \includegraphics[width=13cm]{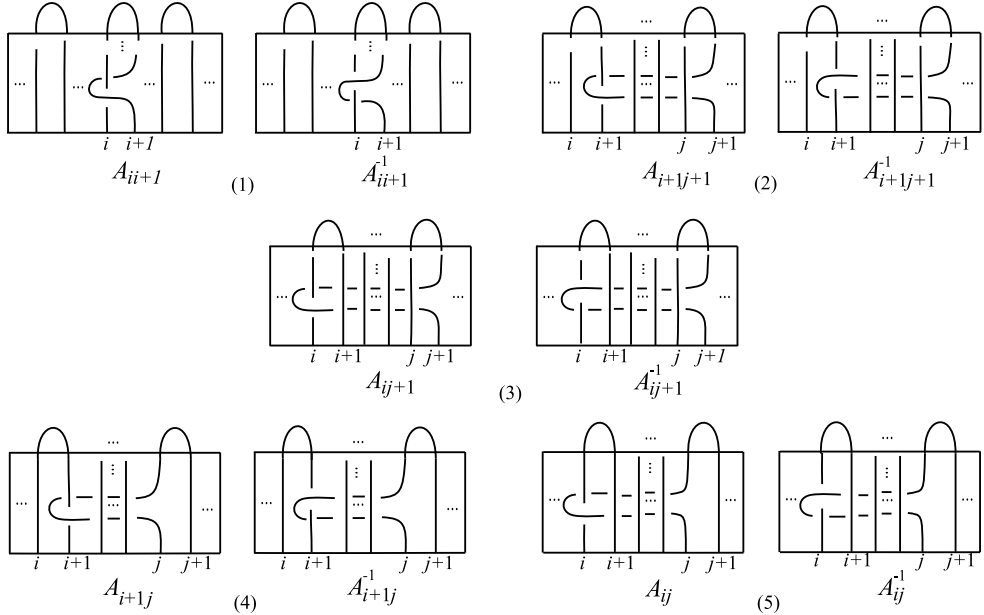}
 \caption{Generators: $A_{ii+1}$, $A_{i+1j+1}$, $A_{ij+1}$, $A_{i+1j}$, $A_{ij}$, $1 \leq i < j \leq 2n$, $i$ and $j$ are both odd.}
  \label{gs}
\end{center}
\end{figure}

\begin{figure}[h!]
\begin{center}
  \includegraphics[width=14cm]{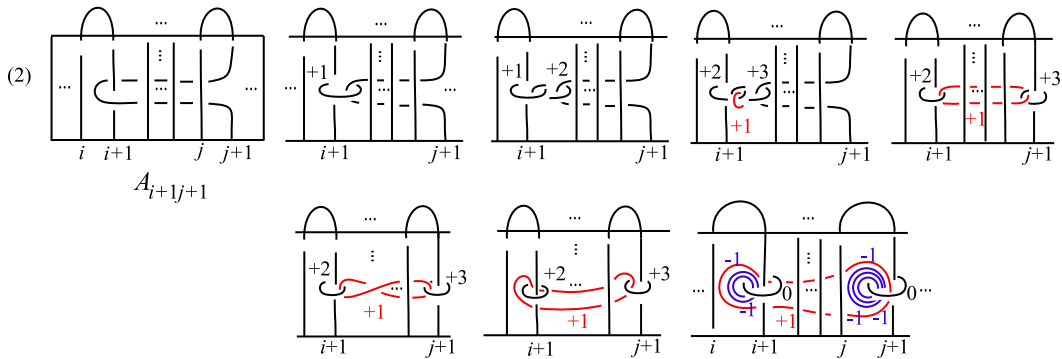}
 \caption{Case (2) $A_{i+1j+1}$.}
  \label{g2}
\end{center}
\end{figure}

\begin{figure}[h!]
\begin{center}
  \includegraphics[width=14cm]{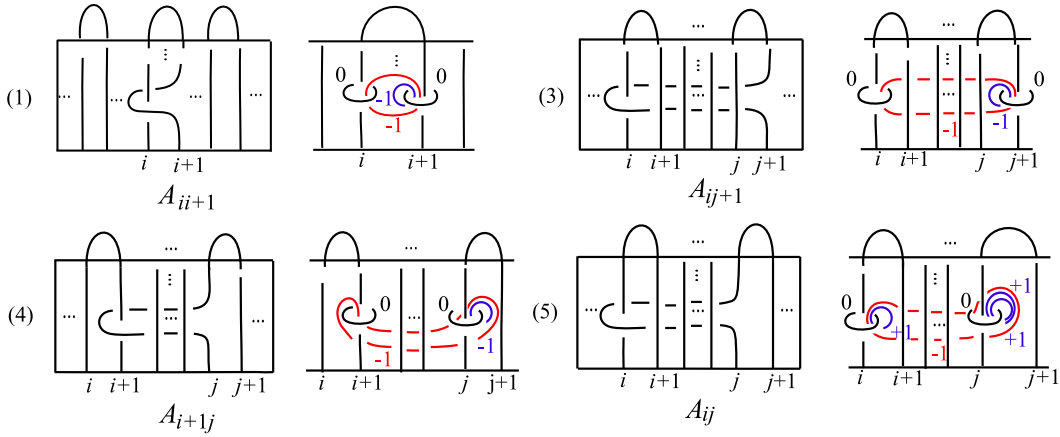}
 \caption{Generators: $A_{ii+1}$, $A_{ij+1}$, $A_{i+1j}$, $A_{ij}$, $1 \leq i < j \leq 2n$, $i$ and $j$ are both odd.}
  \label{gnew}
\end{center}
\end{figure}

 Note that it is enough to verify we can do this for the set of generators and their inverses given in Figure ~ \ref{gs}. All the generators fall in one of the five cases given in Figure ~ \ref{gs}. We explain one complicated case in Figure ~ \ref{g2} and we give a summary for remaining cases in Figure ~ \ref{gnew}. We want to remark that the cases for the inverses are very similar to these cases except the inverses possibly have different framing coefficients. We also want to remark that a pure braided plat presentation of the knot $K$ of the type in Figure ~ \ref{plat}$(b)$ allows us isotope $\pm 1$-framed curves onto a page. For details on each cases, see \cite{thesis}. \begin{flushright}
$\Box$ \end{flushright}

 \begin{theorem} Let $L$ be a link of $k$ components $L_1, \ldots, L_k$ in $S^3$ then $L$ is planar that is $L$ sits on a page of a planar open book for $S^3$.
 \label{top2}
 \end{theorem}
 \begin{proof} Here, we mimic the proof of the Theorem ~ \ref{top}. The only modification required is at the end. Using a pure braided plat presentation of the link $L$, repeatedly blow up to unknot the given link $L$ and arrange the framing of the unknots linking $L$ only once to be $0$ and remove each linking between the unknots linking $L$ to get the middle $\pm 1$-framed curves. After performing the $0$-surgeries, the page of the open book can be constructed by taking the connected sum of components $L_1, \ldots, L_k$ of the link $L$ as shown in Figure ~ \ref{csum}. Hence, we can isotope the middle $\pm 1$-framed curves onto a page using the bands connecting the components. Clearly, the link $L$ sits on a page of this planar open book.
  \begin{figure}[hbt]
 \begin{center}
  \includegraphics[width=10cm]{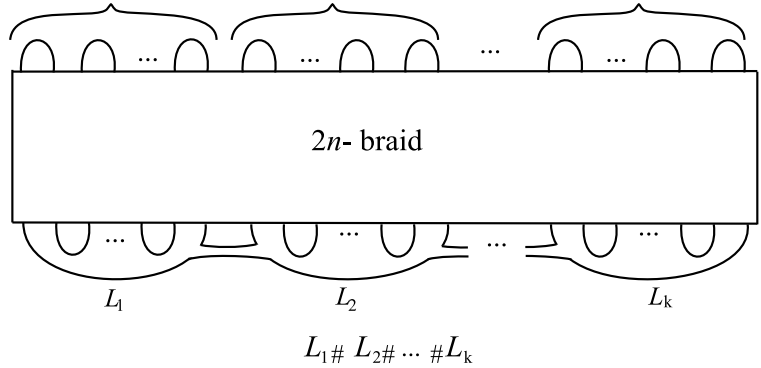}
 \caption{Construct the page of the open book by taking connected sum of the components $L_1, \ldots, L_k$ of the link $L$.}
  \label{csum}
\end{center}
 \end{figure}
 \end{proof}
  We are now ready for the proof of one of the main theorems of this section.  
 \begin{theorem} Let $L$ be a link of $k$ components $L_1, \ldots, L_k$ in a $3$-manifold $M$ then $L$ is planar.
 \label{top3}
 \end{theorem}
 \begin{proof} It is known, see \cite{L} and \cite{W}, that any closed oriented $3$-manifold $M$ may be obtained by $\pm1$ surgery on a link $L_M$ of unknots in $S^3$. Given a link $L$ of $k$ components $L_1, \ldots, L_k$ in a $3$-manifold $M$, we may think of $L$ as a link in $S^3$ which is disjoint from the surgery link $L_M$. Now using the algorithm described in Theorem ~ \ref{top2} we can find a planar open book decomposition for $S^3$ such that the link $L \sqcup L_M$ sits on its page. Also, using a similar idea in Lemma~\ref{rkstab} we can arrange framing of each component of $L_M$ sitting on a page to be  $\pm1$ with respect to the page framing by first stabilizing the open book and then pushing the knot $L$ over the $1$-handle that we use to stabilize the open book. Then away from the link $L$, we can perform $\pm1$ surgeries on $L_M$ which will yield a planar open book for the $3$-manifold $M$ containing the link $L$ on its page. Moreover, this new open book has monodromy which is the old monodromy composed with negative/ positive Dehn twists along each $\pm1$-framed component of the link $L_M$.
 \end{proof}
\section{Support genus of Legendrian Knots}
\begin {definition} The \emph{support genus} $sg(L)$ of a Legendrian knot $L$ in a contact $3$-manifold ($M,\xi$) is the minimal genus of a page of an open book decomposition of $M$ supporting $\xi$ such that $L$ sits on a page of the open book and the framings given by $\xi$ and the page agree.
\end{definition}
\par Given a Legendrian knot $L$ in a contact $3$-manifold $(M, \xi)$, one can always find an open book decomposition compatible with $\xi$ containing $L$ on a page such that the contact framing of $L$ is equal to the framing given by the page. Such an open book decomposition for $(M, \xi)$ can be constructed by an application of Giroux's algorithm, using a contact cell decomposition of $(M, \xi)$ and including the given Legendrian knot $L$ in the $1$-skeleton of the contact cell decomposition, \cite{G}. For Legendrian knots in $(S^3, \xi_{std})$ an alternative algorithm that uses the front projection of Legendrian knots can be found in \cite{AO}, cf. also \cite{F}. Thus the support genus $sg(L)$ of a Legendrian knot $L$ is well defined.
 \begin{theorem} Let $L$ be a null-homologous Legendrian loose knot in an overtwisted contact $3$-manifold ($M,\xi_{ot}$). Then $sg(L) = 0$.
 \label{thmsg}
\end{theorem}
\begin{proof} It is known that two null-homologous Legendrian loose knots $L_1$ and $L_2$ in knot type $K$ with the same Thurston-Bennequin invariant and the same rotation number, then there is a contactomorphism $\psi$ of $(M, \xi_{ot})$ such that $\psi(L_1) = L_2$, \cite{Jknots}. Here, we show that we can realize any pair of integers $(m, n)$ with $m \pm n$ odd as $(tb(L), r(L))$ for a null-homologous loose knot $L$ in knot type $K$ that sits on a planar open book $(S, \varphi)$ supporting ($M,\xi_{ot}$). By Theorem ~ \ref{top3}, we know there is a planar open book decomposition, say $(S_K, \varphi_K)$, for $M$ such that $K$ lies on a page of the open book. The planar open book $(S_K, \varphi_K)$ is compatible with some contact structure $\xi'$ on $M$. If necessary we can negatively stabilize the open book in such a way that the resulting open book is still planar and it is overtwisted. Furthermore, following \cite{Jplanar} we can assume $\xi'$ is the same as the overtwisted contact structure $\xi_{ot}$. Briefly, by performing necessary Lutz twists and taking Murasugi sum of $M$ with an appropriate overtwisted $S^3$, we can arrange the $2$-dimensional invariants $d_2$ and the $3$-dimensional invariants $d_3$ of $\xi'$ and $\xi_{ot}$ to be the same. Thus, the two contact structures will be homotopic, \cite{Gompf}. Then, by Eliashberg \cite{Elias} two overtwisted contact structures will be isotopic. Note that we can do this keeping the open book planar and keeping the given knot $K$ on the page. For the details of how to arrange invariants of overtwisted contact structures, see the proof of Theorem 3.5 in \cite{Jplanar}.
\par Now, we can assume that the planar open book $(S_K, \varphi_K)$ containing the knot $K$ on its page is compatible with the overtwisted contact structure $\xi_{ot}$ on $M$. If necessary by stabilizing the open book positively and pushing the knot $K$ over the $1$-handle, we can assume $K$ is non-separating and we may Legendrian realize $K$ on the page, say it has a Thurston-Bennequin invariant $t'$ and a rotation number $r'$. To realize any pair $(tb(L), r(L))$ for any Legendrian representative of $K$ from the pair $(t', r'),$ first realize the appropriate Thurston-Bennequin invariant $tb(L)$. If $t' > tb(L)$, then to decrease the Thurston-Bennequin invariant stabilize the knot positively or negatively on the page by using Lemma ~ \ref{rkstab}$(1)$. Modify the open book as in Figure ~ \ref{obstab}$(a)$ or $(b)$, both will decrease $tb(L)$. Note, this modification alters neither the contact structure nor the genus of the open book. Now, if $t' < tb(L)$, then to increase the Thurston-Bennequin invariant we need to destabilize the knot positively or negatively on the page by using Lemma ~ \ref{rkstab}$(2)$. Note, this modification alters the contact structure. However, as before, away from the knot by taking the Murasugi sum of $M$ with an appropriate overtwisted $S^3$, we can make sure that the resulting overtwisted contact structure is still isotopic to $\xi_{ot}$.
\par Now, once we realized the pair $(tb(L), r'')$, to complete the proof we only need to realize any possible rotation number $rot(L)$ from $r''$. To increase or decrease the rotation number, we will use Lemma ~ \ref{rkstab} again and stabilize the knot positively or negatively on the page. Recall that a positive and a negative stabilization of a knot increase and decrease the rotation number by $1$, respectively and also recall that both stabilizations decrease the Thurston-Bennequin invariant $tb(L)$ by $1$. Thus, every time we increase or decrease $r''$, we need to make sure that $tb(L)$ stays the same. Clearly, this is possible since to increase the rotation number if we first positively stabilize the knot on the page as in Figure ~ \ref{obstab}$(a)$ and then negatively destabilize the knot on the page as in Figure ~ \ref{obstab}$(d)$, the rotation number will increase by $2$ and $tb(L)$ stays the same. Note that after negatively stabilizing the open book, we again perform a Murasugi sum to keep the contact structure same as $\xi_{ot}$. Similarly, to decrease the rotation number, we first modify the open book as in Figure ~ \ref{obstab}$(b)$ and then as in Figure ~ \ref{obstab}$(c)$, this time the rotation number will decrease by $2$ and $tb(L)$ stays the same. Since $tb(L) \pm rot(L)$ is odd, we can realize any pair $(tb(L), rot(L))$. Thus, for any null-homologous loose Legendrian representative of the knot $K$ we can find a planar open book decomposition supporting $\xi_{ot}$ such that the Legendrian representative sits on the page. 
 \end{proof}

\section{Final Remarks and Questions}

\begin{remark} Note that other than the unknots with $0$-framing coming from resolving the generators $(1)$ $A_{ii+1}$, $(3)$ $A_{ij+1}$, $(4)$ $A_{i+1j}$ in the proof of Theorem ~ \ref{top}, we have only $-1$-framed unknots. In these cases, $-1$-framed unknots contribute positive Dehn twists to the monodromy of the new open book. We want to remark that we can arrange this to be the case for all generators and their inverses. Namely, by blowing up in different ways we can make sure that other than $0$-framed unknots, each case contains only $-1$-framed knots. Thus, at the end we will have an open book decomposition for $S^3$ whose monodromy is a product of only positive Dehn twists and contains the given knot or link on its page. We discuss one case, case $(2)$ $A_{i+1j+1}$, in Figure ~ \ref{rnew}. Other cases can be worked out similarly. For details see \cite{thesis}.
 \begin{figure}[h]
\begin{center}
  \includegraphics[width=15cm]{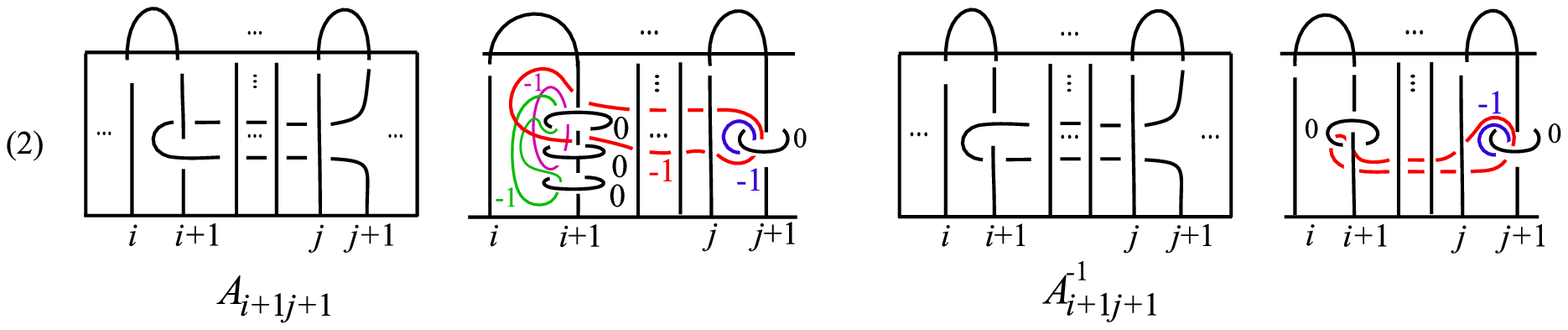}
 \caption{}
  \label{rnew}
\end{center}
\end{figure} 
\end{remark}
As a consequence, we have
\begin{theorem} Any topological knot or link in $S^3$ sits on a planar page of an open book decomposition whose monodromy is a product of positive Dehn twists. 
\end{theorem}
Also, from this observation it follows that
\begin{theorem} Given a knot type $K$ in $(S^3, \xi_{std})$, there is a Legendrian representative $L$ of $K$ such that $sg(L) = 0$.
\end{theorem}
\par It is easy to find examples of support genus non-zero Legendrian knots in weakly fillable tight contact structures.
\begin{lemma} Let $L$ be a Legendrian knot in a weakly fillable tight contact structure with a Thurston-Bennequin invariant $tb(L) > 0$ then $sg(L) > 0$. In particular, any Legendrian knot $L$ in $(S^3, \xi_{std})$ with Thurston-Bennequin invariant $tb(L)\geq 0$ has $sg(L) > 0$.
\label{tbpos}
\end{lemma}
\begin{proof}
\par In \cite{Jplanar}, Etnyre gives constraints on contact structures having support genus zero. In particular, according to \cite{Jplanar} a contact $3$-manifold $(M, \xi)$ obtained by a Legendrian surgery along a Legendrian knot $L$ in a weakly fillable contact structure having Thurston-Bennequin invariant $tb(L) > 0$ has $sg(\xi) > 0$. If a Legendrian knot with $tb(L) > 0$ had support genus $sg(L) = 0$, then performing a Legendrian surgery along $L$ would yield a contact $3$-manifold $(M, \xi)$ with support genus $sg(\xi) = 0$, which is not the case. Therefore, such a Legendrian knot has $sg(L) > 0$. That any Legendrian knot with $tb(L) = 0$ has $sg(L) > 0$ follows from Corollary 1.6 in \cite{OSS}.
\end{proof}
\par There are examples of support genus non-zero non-loose knots in overtwisted contact structures. 
\begin{example} Consider a Legendrian knot $L$ with a Thurston-Bennequin invariant $tb(L) > 0$ in $(S^3, \xi_{std})$. Let $(M, \xi)$ denote the contact $3$-manifold results from a $+1$-contact surgery along a positive stabilization $S_{+}(L)$ of the Legendrian knot $L$. $(M, \xi)$ is overtwisted by \cite{BOZ2} and also c.f by \cite{LS}. Since $tb(S_{+}(L)) \geq 0$ according to the previous Remark ~ \ref{tbpos}, $sg(S_{+}(L)) > 0$. Note that the image $S_{+}(L)'$ of $S_{+}(L)$ in the surgered overtwisted contact manifold $(M, \xi)$ is a non-loose Legendrian knot with a non-zero support genus. The Legendrian knot $S_{+}(L)'$ is non-loose since the complement of $S_{+}(L)'$ in $(M, \xi)$ is contactomorphic to the complement of $S_{+}(L)$ in $(S^3, \xi_{std})$ and $sg(S_{+}(L)') > 0$, otherwise this would contradict to the fact that $sg(S_{+}(L)) > 0$. 
\label{egnonloose}
\end{example}
\par
\begin{figure}[h!]
 \begin{center}
    \includegraphics[width=14cm]{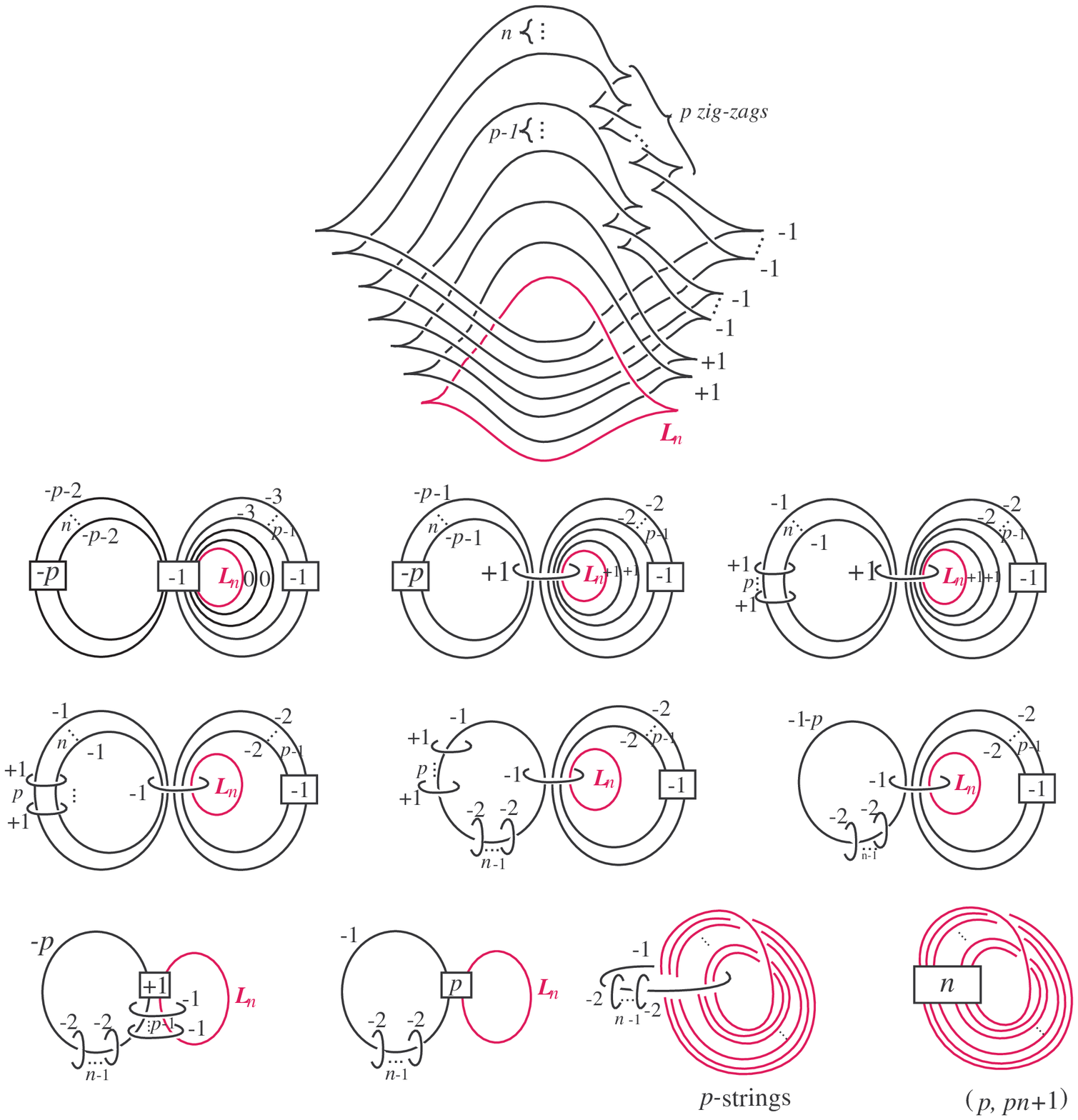}
  \caption{$(p, pn+1)$ Torus knots}
    \label{tk}
 \end{center}
 \end{figure}

 There are examples of support genus zero non-loose knots in overtwisted contact structures. 
\begin{example} The contact $3$-manifold given by the surgery diagram in Figure ~ \ref{tk} is an overtwisted $(S^3, \xi_n)$ with $d_3(\xi_n) = 1 - np(p-1)$, $p > 1$ and $n$ are positive integers. The Legendrian knot $L_n$ in  $(S^3, \xi_n)$ is non-loose with support genus zero and topologically a $(p, pn+1)$ positive torus knot. When $p =2$, Legendrian non-loose knots of knot type $(2, 2n+1)$ positive torus knots first appeared in \cite{LOSS}. Let $X$ denote the $4$-manifold obtained by viewing the integral surgeries as $4$-dimensional $2$-handle attachments to $B^4$. With the help of $X$, we can compute the $3$-dimensional invariant $d_3(\xi_n)$ of the contact structure $\xi_n$. From Figure ~ \ref{tk}, the signature of $X$ is  $\sigma(X) = -n-p+1$ and the Euler characteristic of $X$ is $\chi(X) = n + p + 1$. Also, using a second homology class $c \in H^2(X, \mathbb{Z})$ defined by the rotation number, we compute $c^2 = -n(2p - 1)^2 - (p-1)$. From the formula: 
\begin{equation*}
d_3(\xi) = \frac{1}{4}(c^2 - 3(\sigma(X)) - 2(\chi(X))) + q
\end{equation*}
where $q$ denotes the number of $+1$-contact surgeries, the $3$-dimensional invariant of $\xi_n$ is $d_3(\xi_n) = 1 - np(p-1)$. Note that $\xi_n$ is overtwisted since $d_3(\xi_n) < 0$. Note also that $L_n$ is non-loose since Legendrian surgery on $L_n$ cancels one of the $+1$-surgeries in Figure ~ \ref{tk} and results in a tight contact structure. By a similar argument used in \cite{SB}, the surgery link together with the Legendrian knot $L_n$ given in Figure ~ \ref{tk} can be put on a page of a planar open book of $(S^3, \xi_{std})$. After performing surgeries, we will get $(S^3, \xi_n)$ compatible with a planar open book containing the Legendrian knot $L_n$ on its page. Therefore, $sg(L_n) = 0$. 
\end{example}
 \begin{remark}As we discussed in Example ~ \ref{egnonloose} above, in overtwisted contact structures there are examples of non-loose knots having support genus non-zero. Let $L$ be a null-homologous, support genus non-zero, non-loose Legendrian knot of knot type $K$ in an overtwisted contact manifold $(M, \xi_{ot})$. We can find a loose knot $\tilde{L}$ of knot type $K$ in $(M, \xi_{ot})$ such that $\tilde{L}$ has the same classical invariants as $L$. Moreover, by Theorem ~ \ref{thmsg} it follows that $sg(\tilde{L}) = 0$. Thus, we have examples of knots having the same classical invariants but different support genus in overtwisted contact structures. Also, it would be very interesting to know
\begin{question} What are the examples of Legendrian knots $L_1$ and $L_2$ in a tight contact structure such that $tb(L_1) = tb(L_2)$ and $rot(L_1) = rot(L_2)$ but $sg(L_1) \neq sg(L_2)$?
\end{question}
\end{remark}

As explained in Lemma~ \ref{rkstab}, if a Legendrian knot $L$ sits on a page of an open book decomposition, then positive or negative stabilization of $L$ can be realized on the page of the open book as in Figure ~ \ref{obstab}. Note that we add $1$-handles in such a way that the resulting open book still has the same genus. As a result, we have
\begin{theorem} If a Legendrian knot $L$ has support genus $sg(L) = n$, then the stabilizations $S_{+}^{n_1} S_{-}^{n_2}(L)$ of $L$ have support genus $sg(S_{+}^{n_1}S_{-}^{n_2}(L)) \leq n$.
\label{maxgenus}
\end{theorem} 
\par By the above Theorem ~ \ref{maxgenus}, given a knot type $K$, if all Legendrian knots realizing $K$ with Thurston-Bennequin invariant strictly less than the maximal Thurston-Bennequin invariant destabilize and the Legendrian knots with maximal Thurston-Bennequin invariant has support genus zero, then all Legendrian knots of knot type $K$ has support genus zero. For example, all Legendrian unknots in $(S^3, \xi_{std})$ are planar. 
\par Recall that for a non-zero rational number $r \in \mathbb{Q}$, a contact $r$-surgery on a Legendrian knot $L$ in a contact $3$-manifold $(M, \xi)$ is a topological $r$-surgery with respect to the contact framing. The resulting manifold is a new contact $3$-manifold $(M', \xi')$ where the contact structure $\xi'$ is constructed by extending $\xi$ from the complement of a standard contact neighborhood of $L$ to a tight contact structure on the glued solid torus, \cite{DG}. Such an extension always exists and it is unique when $r = \frac{1}{k}$, $k \in \mathbb{N}$, \cite{Ko}. 
\begin{theorem} Let $L$ be a Legendrian knot in a contact $3$-manifold $(M, \xi)$ such that $sg(L) = 0$. Then, the contact $3$-manifold $(M', \xi')$ results from a contact $r$-surgery along $L$ has $sg(\xi') = 0$. 
\end{theorem}
\par We want to remark that rational contact surgeries on a Legendrian link in $(S^3, \xi_{std})$ on pages of open book decompositions first discussed in \cite{BOz}.
\begin{proof} For contact $r$-surgery with $r < 0$, consider a continued fraction expansion of $r - 1$ 
\begin{equation*} 
[r_1, r_2, \ldots, r_n] = r_1 - \frac{1}{r_2 - \frac{1}{\Large{\cdots - \frac{1}{r_n}}}}
\end{equation*}
with integers $r_i \leq -2$, $i = 1, \ldots, n$. Let $L_1$ be the $|r_1 + 1|$ times stabilization of the front projection of the Legendrian knot $L$ and let $L_i$ be the Legendrian push off of $L_{i-1}$ with additional $|r_i + 2|$ stabilizations, $i = 2, \ldots, n$. Then following \cite{DG}, we can replace contact $r$-surgery along $L$ by a sequence of contact $-1$-surgeries along $L_1, \ldots, L_n$. Since the support genus $sg(L) = 0$, by Lemma~\ref{rkstab}$(1)$ and by keeping the page of the open book planar we can realize each $L_i$ on a planar open book containing $L$ on its page. Again by using Lemma~\ref{rkstab}$(1)$ we can arrange framing of each $L_i$ sitting on a planar page to be $-1$ with respect to the page framing. After performing contact surgeries, we will obtain a support genus zero contact $3$-manifold.
\par According to \cite{DG}, for $p$, $q$ relatively prime positive integers, a contact $r = \frac{p}{q} > 0$ surgery on $L$ corresponds to $k$ contact $+1$-surgeries along $k$ Legendrian push offs of $L$ followed by a contact $r' = \frac{p}{q-kp}$-surgery along a Legendrian push off of $L$ for any integer $k \in \mathbb{N}$ such that $q - kp < 0$. By starting with a planar open book containing the Legendrian knot $L$ on its page, we can easily see $k$ Legendrian push offs of $L$ on the page. By using Lemma~\ref{rkstab}$(1)$, we can arrange the framings of each push off of $L$ sitting on a planar page to be $+1$ with respect to the page framing. Hence to complete the proof we need only to show we can perform $r' < 0$ surgery on a Legendrian push off of $L$ on the page also but this can be easily arranged as we did above. 
\end{proof}

\begin{remark}Note that the support genus of a Legendrian knot gives an upper bound on the support genus of a contact structure, that is, $sg(L) \geq sg (\xi)$. So, if there is a Legendrian knot $L$ in a contact $3$-manifold $(M, \xi)$ having support genus zero, then $sg(\xi) = 0$. Also, it is still not known if there is a contact structure $\xi$ having support genus $sg(\xi) > 1$. This raises a natural question:
\begin{question} Is $sg(L)$ or $sg(\xi)$ ever bigger than $1$? 
\end{question}
\end{remark}
\begin{remark} In \cite{LOSS}, Heegaard Floer invariants of Legendrian knots are defined by using an open book decomposition that supports the contact structure and contains the knot on its page. Also, in \cite{OSS}, Heegaard Floer homology is used to give restrictions on contact structures supported by planar open book decompositions. One can investigate the work on knots sitting on planar pages in this paper from Heegaard Floer theory perspective:
\begin{question} Can one find restrictions on Legendrian knots sitting on planar pages by using Heegaard Floer homology? 
\end{question}
\end{remark}
Here are some final questions:
\begin{question} Let $L$ be a Legendrian knot in $(S^3, \xi_{std})$ with $tb(L) < 0$. Is $sg(L) = 0$?
\end{question}
\begin{question} Is the support genus of knots additive under connected sums?
\end{question}
 \begin{question} What is the relation between support genus of a knot and its mirror? 
\end{question} 

\section{Appendix}
In this Appendix, we prove Lemma~\ref{rkstab} in Section 2. Throughout we use the standard terminology from convex surface theory. See \cite{Ko} and \cite{EH3} for detailed background.
\par By looking at the characteristic foliation on a disk cobounded by a Legendrian knot and a stabilization of the Legendrian knot we may see how to destabilize the Legendrian knot.
\par Recall that the \textit{characteristic foliation} $\Sigma_{\xi}$ of a surface $\Sigma$ is the singular foliation induced on $\Sigma$ from $\xi$ where $\Sigma_{\xi}(p) = \xi_p \cap T\Sigma_p $, $p \in \Sigma$. The singular points are the points where $\xi_p = T\Sigma_p$. Any surface $\Sigma$ may be perturbed so that its characteristic foliation $\Sigma_{\xi}$ has only generic isolated singularities, elliptic singularities and hyperbolic singularities. The singularity is \textit{positive} (\textit{negative}) if the orientation on $\xi_p$ agrees (disagrees) with the orientation of $T\Sigma_p$. 
\par Also recall that a closed oriented surface $\Sigma$ in a contact manifold $(M,\;\xi)$ is called \textit{convex} if there is a \textit{contact vector field v}, that is a vector field whose flow preserves the contact structure $\xi$, transverse to $\Sigma$. Given a convex surface $\Sigma$ in (M,\;$\xi$) with a contact vector field $v$, the \textit{dividing set} $\Gamma_{\Sigma}$ of $\Sigma$ is defined as $\Gamma_{\Sigma} =\{ x \in \Sigma : v(x) \in \xi_x\}$.

\par
\begin{figure}[h!]
\begin{center}
  \includegraphics[width=7cm]{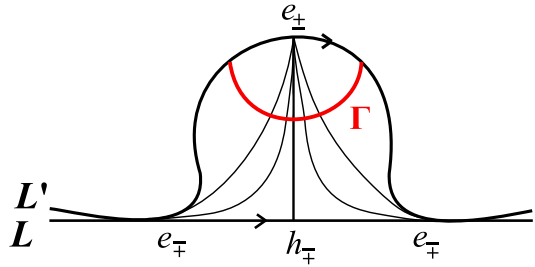}
 \caption{A bypass for $L$ and $L' = S_{\pm}(L)$. The curve $\Gamma$ is the dividing curve of $D$, $e_{\pm}$ denotes positive/ negative elliptic singularity and $h_{\pm}$ denotes positive/ negative hyperbolic singularity.}
  \label{destab}
\end{center}
\end{figure}

 Given an oriented Legendrian knot $L$ the positive stabilization $S_{+}(L)$ of $L$ and the Legendrian knot $L$ cobound a convex disk $D$ where $tb(\partial D) = -1$ and $D \cap L$ contains two negative elliptic and one negative hyperbolic singularities and $D \cap S_{+}(L)$ contains the same two negative elliptic singularities and one positive elliptic singularity. Similarly, the negative stabilization $S_{-}(L)$ of $L$ and the Legendrian knot $L$ cobound a convex disk $D$ where $tb(\partial D) = -1$ and $D \cap L$ contains two positive elliptic and one positive hyperbolic singularities and $D \cap S_{-}(L)$ contains the same two positive elliptic singularities and one negative elliptic singularity. Such a disk is called a \textit{stabilizing disk} for $L$ or a \textit{bypass} for $L$ and $S_{\pm}(L)$. See Figure~\ref{destab}. Note that all the singularities of $D_{\xi}$ have the same sign except one which indicate us whether we are positively or negatively stabilizing the Legendrian knot $L$. For a detailed discussion of stabilizations and bypass disks see \cite{J2}, \cite{EH3}. 
\par {\em \textbf{Proof of Lemma ~ \ref{rkstab}}.} $\textbf{(1)}$ To prove $(1)$ we find a stabilizing disk for each case as we discussed above. See Figure~\ref{destab}. First, positively stabilize the open book as in Figure~\ref{obstab}$(a)$ and push the Legendrian knot $L$ over the $1$-handle that is used to stabilize the open book positively, call the new curve $L_{+}$. We will show that $L_{+}$ is a positive stabilization $S_{+}(L)$ of $L$. 
\par
\begin{figure}[h!]
\begin{center}
  \includegraphics[width=13cm]{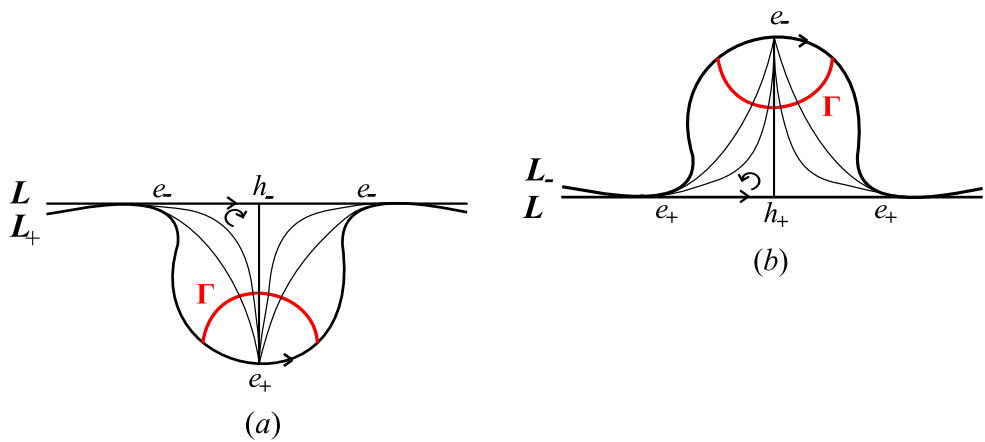}
 \caption{$(a)$ Positive stabilization $L_{+} = S_{+}(L)$ of $L$, $(b)$ Negative stabilization $L_{-} = S_{-}(L)$ of $L$.}
  \label{obstab1}
\end{center}
\end{figure}

 Notice the Legendrian unknot $a$ with $tb(a) = -1$ in Figure~\ref{obstab}$(a)$. Legendrian unknot $a$ bounds a disk $D$ in $M$. Since $tb(a) = -1$, $D$ is convex and the dividing curves intersect $\partial D$ twice. Now, we can think $L_{+}$ as the knot obtained from pushing $L$ across $D$. Note that $D$ is a bypass for $L$ and $L_{+}$. See Figure~\ref{obstab1}$(a)$, the curve $\Gamma$ denotes the diving curve of $D$. A singularity along $\partial D$ is positive or negative depending on whether the contact planes passing $D$ are twisting in a right handed fashion or a left handed fashion. The sign of the singularities is determined by using the orientation of $L$ which determines the orientation of $D$ near the boundary. See Figure~\ref{obstab1}$(a)$ again.
\par The negative stabilization $S_{-}(L)$ of the Legendrian knot $L$ can be realized on a page of the open book decomposition in a similar way. This time we use the Legendrian unknot $b$ with $tb(b) = -1$ in Figure~\ref{obstab}$(b)$ and the convex disk that $b$ bounds in $M$. See Figure~\ref{obstab1}$(b)$.
\bigskip
\par $\textbf{(2)}$ We prove $(2)$ for null-homologous Legendrian knots only. First, negatively stabilize the open book as in Figure~\ref{obstab}$(c)$ and then push the knot $L$ over the $1$-handle that is used to stabilize the open book negatively, call the new curve $L'_{+}$. 
\par Note that in general the negative stabilization of an open book decomposition changes the contact structure $\xi$. However, in this case the curve $L$ on the page gives a Legendrian knot $L'$ in the new contact structure. We will show that $L'_{+}$ in Figure~\ref{obstab}$(c)$ is a positive destabilization of $L'$. 
\par
\begin{figure}[h!]
\begin{center}
  \includegraphics[width=15cm]{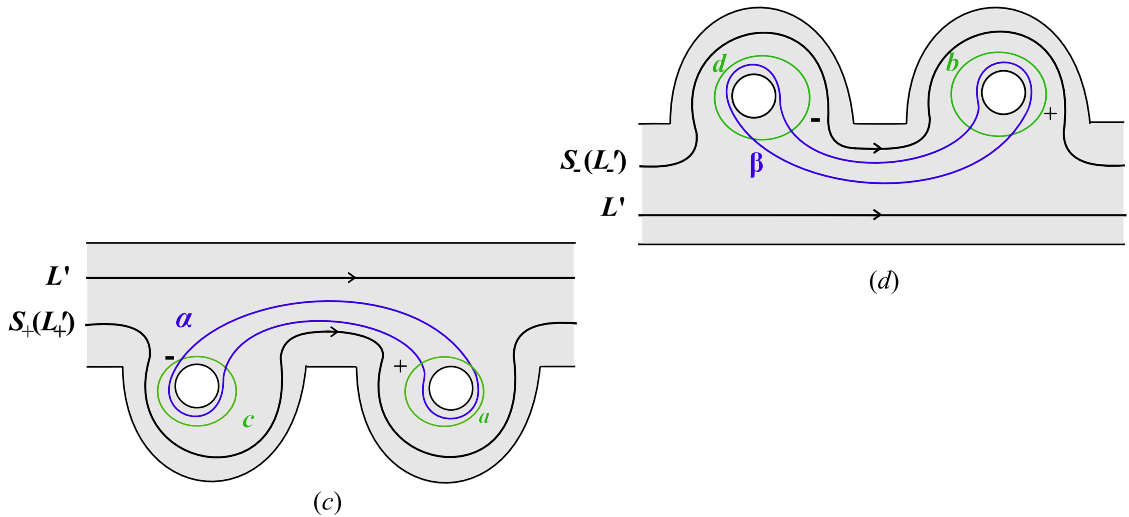}
 \caption{$(c)$ Positive destabilization of $L'$, $S_{+}(L'_{+}) = L'$, $(d)$ Negative destabilization of $L'$, $S_{-}(L'_{-})= L'$.}
  \label{obstab2}
\end{center}
\end{figure}

 We want to remark that the Legendrian unknot $c$ with $tb(c) = +1$ in Figure~\ref{obstab}$(c)$ bounds a disk $D$ in $M$. Since $D$ is not convex unlike in the proof of $(1)$ we can not use this disk to find a bypass. Instead, we positively stabilize the open book as in Figure~\ref{obstab2}$(c)$ and push the Legendrian knot $L'_{+}$ over the $1$-handle that we use to stabilize the open book positively. By $(1)$, the resulting Legendrian knot is a positive stabilization $S_{+}(L'_{+})$ of $L'_{+}$. We will show that $S_{+}(L'_{+})$ is Legendrian isotopic to $L'$. Note that the curve $\a$ in Figure~\ref{obstab2}$(c)$ is a Legendrian unknot with $tb(\a) = 0$. In fact, Legendrian unknot $\a$ bounds an overtwisted disk which is disjoint from $L'_+$ in $M$. Legendrian knots $L'$ and $S_{+}(L'_{+})$ have the same classical invariants, that is, they have the same knot type, same Thurston-Bennequin invariant and same rotation number, and since they have a common overtwisted disk in their complement, by \cite{Dymara} $L'$ and $S_{+}(L'_{+})$ are Legendrian isotopic.
\par Similarly, the negative destabilization $L'_{-}$ of the Legendrian knot $L'$ can be realized on a page of the open book decomposition. We stabilize the open book as in Figure~\ref{obstab2}$(d)$ and push the Legendrian knot $L'_{-}$ over the $1$-handle to get negative stabilization $S_{-}(L'_{-})$ of  $L'_{-}$. We conclude that $S_{-}(L'_{-})$ and $L'$ are Legendrian isotopic by using the Legendrian unknot $\beta$ which has $tb(\beta) = 0$ in Figure~\ref{obstab2}$(d)$. \begin{flushright} $\Box$ \end{flushright}

\par \textbf{Acknowledgments.} I am grateful to John B. Etnyre for his invaluable support and guidance. I would like to thank Stavros Garoufalidis, \c{C}a\u{g}r{\i} Karakurt and Kenneth L. Baker for helpful conversations. The author is partially supported by NSF grant DMS-0804820 and by TUBITAK-2214.

\end{document}